\documentclass[12pt]{article}

\usepackage[utf8]{inputenc}

\usepackage{amsfonts}
\usepackage{amssymb}
\usepackage{amsmath}
\usepackage{amsthm}
\usepackage[T2A]{fontenc}

\usepackage{amsfonts}
\usepackage{amssymb}
\usepackage{amsmath}
\usepackage{amsthm}
\usepackage{color}
\begin{document} 

\begin{Huge}

    \centerline{Liminf-results for sums with Kronecker sequence}
    \end{Huge}
    \begin{Large}
    \vskip+0.5cm
    \centerline{A.O. Chebotarenko}
\end{Large}
\vskip+2cm

Let $\theta \in \mathbb{R}$ be irrational number and $f:\mathbb{R}\to \mathbb{R}$ be 1-periodic function with zero mean value
\begin{equation}\label{zeromen}
\int_0^1 f(x)dx=0.
\end{equation}
In the present paper we study the behavior of the sums 
 \begin{equation}\label{sum}
\sum_{k=0}^{Q-1}
f(k\theta+\varphi) 
 \end{equation}
related to  the Kronecker sequence 
 \begin{equation}\label{kroo}
 \{k\theta\},\,\,\ k =0,1,2,3,... \, .
 \end{equation}
\section{Introduction and main results}

 Here in we introduce the main object of our research and formulate our results.

 Let  $f(x)$ be a 1-periodic function of bounded variation 
 and $\xi_0,....,\xi_{Q-1} \in [0,1]$. The famous  Koksma's inequality (see \cite{KN}, Ch.2, \S 5)  ensures the bound
 \begin{equation}\label{ko}
 \left|\sum_{k=0}^{Q-1} f(\xi_k) -  Q \cdot \int_0^1f(x)dx\right|
 \le {\rm Var}[f] \cdot  D (\Xi),
 \end{equation}
 where ${\rm Var}[f]$ is the variation of function $f$
 and 
 $$
 D(\Xi) =
 \sup_{\gamma\in [0,1]}
 \left|
 \# \{ j\in\{ 0,...,Q-1\}: \xi_j \le \gamma\} - Q\gamma
 \right| 
 $$
 is the discrepancy of the set $ \Xi = \{ \xi_0,...,\xi_{Q-1}\}$.

\vskip+0.3cm
   For the  analysis of uniform distribution of the sequence (\ref{kroo})   the approximations of $\theta$  by its {\it convergent fractions} are of major importance.
 Recall that every irrational $\theta$ may be uniquely represented as ordinary continued fraction
 $$
 \theta = 
 [a_0;a_1,a_2,...,a_\nu,...]
   =   a_0  +
\frac{1}{\displaystyle{a_1+\frac{1}{\displaystyle{a_2 + \cdots+
\frac{1}{\displaystyle{a_\nu + ...
{} }}}}}}
   , a_0\in \mathbb{Z},\, a_\nu\in \mathbb{Z}_+, n=1,2,3,... \,\,,
 $$
 and convergent fractions for $\theta$ are defined as
 $$
 \frac{P_\nu}{Q_\nu} =    [a_0;a_1,a_2,...,a_\nu].
 $$
 Let $D_Q$ be the
  discrepancy of the first $Q$ elements of the Kronecker sequence (\ref{kroo}).
  When $ Q = Q_\nu$ is a denominator of a convergent to $\theta$, one has
 $ D_{Q_{\nu}} =O(1).$ So for  function $f$ of bounded variation and irrational $\theta$ by (\ref{ko}) we have
 \begin{equation}\label{ooBV}
 \sup_\nu  \max_{\varphi \in [0,1]} \left|\sum_{k=0}^{Q_\nu-1} f(k\theta+\varphi)\right|<+\infty.
 \end{equation}
 The study of liminf-results for the absolute values of the sums (\ref{sum}) was initiated by Kozlov \cite{koz}.
The following statement  proved by Sidorov  \cite{sy} gives a more precise result  than (\ref{ooBV}) in the case when function $f(x)$ is 
{\it absolutely continuous} (this condition is stronger then the condition of bounded variation, see \cite{KF}, Ch. 9).
 \vskip+0.3cm
 
 {\bf Theorem A.}
{\it 
 For any  1-periodic absolutely continuous function $f(x)$  with  zero mean value and any irrational $\theta$ one has
  $$
 \liminf_{Q\to \infty} \max_{\varphi \in [0,1]} \left|\sum_{k=0}^{Q-1} f(k\theta+\varphi)\right| =
 \lim_{\nu\to \infty} \max_{\varphi \in [0,1]} \left|\sum_{k=0}^{Q_\nu-1}
f(k\theta+\varphi)\right|=0.
 $$
}

A brief survey of liminf-results for the  sums (\ref{sum}) as well as some multidimensional generalizations one can find in a recent paper \cite{mo}.
 
 \vskip+0.3cm
 
This paper is motivated by a question by R. Tichy \cite{T}, who asked if a converse statement to Theorem A is true.
In the present paper we answer the question by Tichy and prove several new results concerning sums (\ref{sum}).
Our first theorem  shows that $\liminf$-result may hold of all irrational $\theta $ and for all $\varphi$
even for function $f$ which is not  a function of bounded variation.
  \vskip+0.3cm
{\bf Theorem 1.}
{\it
There exists a continuous  1-periodic function $f(x)$ which is not a function of bounded variation such that  for any  irrational $\theta\in \mathbb{R}$ and for any $\varphi \in \mathbb{R}$  the asymptotic equality
 \begin{equation}\label{0}
 \liminf_{Q\to \infty} \left| \sum_{k=0}^{Q-1} f(k\theta +\varphi) \right| = 0
 \end{equation}
 is valid.}
 \vskip+0.3cm

  \vskip+0.3cm
 {\bf Remark 1.} 
 Of course from the uniform distribution property of the sequence (\ref{kroo})  it follows that the function from Theorem 1 should satisfy zero mean value condition
 (\ref{zeromen}).

   \vskip+0.3cm
 
 Our second result shows that the conclusion of Theorem A is not valid under a weaker condition that 
 continuous $f$ has bounded variation. 
 The author does not know any such example documented in literature before.

\vskip+0.3cm
 
{\bf Theorem 2.}
{\it Let irrational number $\theta = [0;,1,1,a_3,a_4,...]$  satisfy the condition
\begin{equation}\label{oo1}
\lim_{\nu\to \infty}\frac{Q_{\nu-1}}{a_{\nu+1} }
=0.
\end{equation}
Then there exists  $1$-periodic continuous function $g_\theta$ of bounded variation and with zero mean value
\begin{equation}\label{oozero}
\int_0^1 g_\theta (x)dx = 0
\end{equation}
such that   
\begin{equation}\label{oozero1}
\liminf_{\nu\to\infty}\sum_{k=0}^{Q_\nu-1}g_\theta(k\theta) > 0.
\end{equation}
}
 \vskip+0.3cm

Theorem  2 states that for any irrational number which has very good rational approximations
(for example, to satisfy (\ref{oo1}) one can take $\theta$ such  that for approximation by convergents one has
$$
\left|\theta - \frac{P_\nu}{Q_\nu}\right| \le \frac{1}{Q_\nu^{3+\varepsilon}}, \,\,\, \varepsilon >0
$$
for every $\nu $)
we can construct a continuous function of bounded variation with property (\ref{oozero1}). In addition to such a result we formulate below a  simple statement in converse direction which goes back to some constructions from \cite{koz}. Recall  that {\it modulus of continuity} for function $f:\mathbb{R}\to \mathbb{R}$ is defined as
$$\omega(\delta) = \sup\{|f(x) - f(y)|: |x-y|<\delta, \ x,y \in \mathbb{R}\}.$$
Function $\delta\mapsto \omega(\delta)$ is  an increasing function
and for    uniformly continuous  function $f$
 we have
 $$\lim_{\delta \to 0} \omega(\delta) = 0.$$
Here we deal with continuous 1-periodic functions and such functions, of course, are uniformly continuous.

 \vskip+0.3cm
 {\bf Proposition 1.}
 {\it 
Let $f$ be continuous $1$-periodic function of bounded variation with zero mean value (\ref{zeromen}). 
Suppose that irrational $\theta$ satisfies  the condition
\begin{equation}\label{oo1o}
\liminf_{\nu\to \infty} Q_\nu \cdot 
\omega (||Q_\nu \theta||)  =0.
\end{equation}
Then   for any $\varphi \in \mathbb{R}$ one has 
\begin{equation}\label{ooo1}
    \liminf_{N\to \infty} \,  \sum_{k=0}^{N-1} f(k\theta+\varphi)  = 0.
\end{equation}

}

 \vskip+0.3cm
 {\bf Remark 2.} If we consider  periodic continuous function $f(x)$, condition
 (\ref{oo1o}) means that $\theta$ has partial quotients which grow fast enough. In particular, for Lipschitz $f$, that is in the case when $ \omega (\delta) \leq c\delta$ with some positive $c$, condition  (\ref{oo1o}) means  that partial quotients of $\theta$
 are just unbounded.

 \vskip+0.3cm

 {\bf Remark 3.} In  Proposition 1 we do not claim that asymptotic equation (\ref{ooo1})
holds uniformly in $\varphi$.
  \vskip+0.3cm

Finally, we formulate one more simple statement dealing with {\it even} functions.

  \vskip+0.3cm
 {\bf Proposition 2.}
 {\it 
 Let $f$ be continuous $1$-periodic function of bounded variation with mean value (\ref{zeromen})  and
\begin{equation}\label{even}
f(x) = f(-x) \  \forall x \in \mathbb R.
\end{equation}
Then for any irrational $\theta$  one has
\begin{equation}
    \liminf_{Q\to \infty} \left| \sum_{k=0}^{Q-1} f(k\theta) \right| = 0.
\end{equation}}

  \vskip+0.3cm
 The structure of the rest of the paper is very simple. In Section 2 we prove Theorem 1, in Section 3 we prove Theorem 3 and in Sections 4,5 we prove Propositions 1,2 respectively.

\section{Proof of Theorem 1}

In this section we give a proof of Theorem 1.
First of all in Subsection \ref{f} we introduce the construction of function $f(x)$
which in some sense generalizes
the famous example $ x\mapsto x\sin \left(\frac{1}{x}\right)$ of a function of unbounded variation.
In Subsection \ref{f3} we study our sums by means of approximation $\theta$ by convergents.
In Subsection \ref{f4} we finalize the proof of Theorem 1.

\subsection{Construction of function $f$}\label{f}

We define numbers $ \varepsilon_n, x_n$ and  absolutely continuous functions $ h_n(x), f_n (x)$, $ n =1,2,3,...$ 
with zero mean value
$
\int_0^1 f(t) dt = 0
$
by induction.

First of all we define 
$$
x_0 = 0, \,\,\,\, x_1 = \varepsilon_1 = \frac{1}{2}
$$
and 1-periodic absolutely continuous functions 
$$
f_1(x) = h_1(x) = \sin ( 4\pi x) \cdot {\bf 1}_{[x_0,x_1]} (x),\,\,\,\, 0\le x \le 1.
$$
Assume that  $ \varepsilon_j, x_j$ and functions $ h_j(x), f_j (x)$, $ j =1,2,3,..., n$ are defined.
Then we consider  $f_n(x)$. It is an absolutely continuous function.  So for its derivative $f_n'(x)$ we have $
f_n'(x) \in L_1$. As the set of trigonometric polynomials is dense in $L_1$ (see \cite{KF}),    there exists a trigonometric polynomial
$$
p_n(x) = \sum_{m\in \mathbb{Z}: \, 0<|m|\le N_n} p_{n,m} e^{2\pi i mx},\,\,\,\,\, N_n \ge 1
,$$
such that 
\begin{equation}\label{1}
\int_0^1 | f_n'(t) - p_n'(t)|dt \le \frac{1}{n}.
\end{equation}
Now we define
$$
M_{-1}=M_0 = 1,\,\,\,
M_n = M[p_n] = \sum_{m\in \mathbb{Z}: \, 0<|m|\le N_n} |mp_{n,m}|
$$
and put
\begin{equation}\label{2}
\varepsilon_{n+1} =  \varepsilon_n \, \min \left( \frac{1}{2N_n}, \frac{1}{n M_{n-2}}\right) ,\,\,\,\,
x_{n+1} = x_n +\varepsilon_n,\,\,\, n = 1,2,3,... 
\end{equation}
and
$$
h_{n+1}(x) = \frac{1}{n+1} \sin \left(2 \pi \left(\frac{x-x_n}{\varepsilon_{n+1}}\right)\right) \cdot
 {\bf 1}_{[x_n,x_{n+1}]} (x),\,\,\,\,
 f_{n+1}(x) = f_n(x) + h_{n+1}(x).
$$
It is clear from the definition (\ref{2}) that 
$\varepsilon_{n+1}\le \frac{\varepsilon_n}{2}$ and so
\begin{equation}\label{ee}
\sum_{j=n+1}^\infty \varepsilon_j \le \varepsilon_n
\end{equation}
and
\begin{equation}\label{xx}
x_n = \sum_{j=1}^n \varepsilon_j < 1,\,\,\,\,\,  x^* = \sum_{j=1}^\infty \varepsilon_j \le 1
\end{equation}
Notice that the supports of $h_n$ is just the interval $(x_{n-1}<x_n)$ and for every $n$ we have
\begin{equation}\label{var}
{\rm Var} [h_n] = \frac{4}{n}
.
\end{equation}
Now we define
$$
f(x) = \sum_{j=1}^\infty h_j(x).
$$
It is clear 
that
$$
{\rm Var} [f]
=
\sum_{j=1}^\infty
{\rm Var} [h_j]
=+\infty.
$$
So $f(x)$ is a continuous function with zero mean value which is not a function of bounded variation.

\subsection{Approximation by convergents and auxiliary statements}\label{f3}

For irrational $\theta$ we consider its continued fraction expansion 
$\theta = [a_0;a_1,a_2,..., a_\nu],...$ and convergent fractions
$\frac{P_\nu}{Q_\nu} =  [a_0;a_1,a_2,..., a_\nu]$. 
We need two well known facts (see \cite{CF}). The first is that 
\begin{equation}\label{12}
Q_{\nu+1} ||Q_\nu \theta|| >
\frac{1}{2}\,\,\,\, \forall \nu.
\end{equation}
The second is that the denominators $Q_\nu$ form a sequence of all the best approximations,
that is
\begin{equation}\label{be}
||Q_\nu \theta|| < || q\theta||\,\,\,\,\,\,\, \text{for all  }\,\,\, q\in \mathbb{Z}_+\,\,\,
\text{such that }\,\,\,
1\le |q| < Q_{\nu+1}, \,\,q\neq \pm Q_\nu
\end{equation}

Here we prove two easy lemmas.

\vskip+0.3cm
{\bf Lemma 1.}
{\it  Let  $I\subset [0,1)$ be an interval and $g(x)$ be 1-periodic function such that $g(x) = 0$ for all $x\not\in I$.
Assume that $ ||Q_{l-1}\theta|| > \frac{|I|}{2}$. Then for any $\varphi$ one has 
$$
\left|
\sum_{k=0}^{Q_l-1} g(k\theta+\varphi)\right|
\le 2 \sup_{x\in [0,1]} |g(x)|.
$$}

\begin{proof}
Assume that $ 0\le k'< k''\le Q_l-1$. Then $ 0< k''-k'< Q_l$ and as denominators of continued fractions give best approximations we deduce  from (\ref{be}) that 
$$ 
|(k''\theta+\varphi)-
(k'\theta+\varphi)|\ge ||(k''-k')\theta|| \ge ||Q_{l-1}\theta|| > \frac{|I|}{2}.
$$
So in the interval $I$ there exist not more than two  points of the form $\{k\theta+\varphi\}$.
\end{proof}

 {\bf Lemma 2.}
{\it  Let 
$$
p(x) = \sum_{m\in \mathbb{Z}: \, 0<|m|\le N} p_m e^{2\pi i mx}
$$
be a trigonometric polynomial,
$$ 
M[p] =   \sum_{m\in \mathbb{Z}: \, 0<|m|\le N} |m p_m|,
$$
and $ N < Q_{\nu+1}$.
Then for any $\varphi$ one has 
$$
\left|
\sum_{k=0}^{Q-1} 
p(k\theta +\varphi) 
\right|
\le
\frac{\pi}{2}\, 
M[p] \cdot \frac{||Q\theta||}
{||Q_{\nu}\theta||}.
$$}
\begin{proof}
$$
\left|
\sum_{k=0}^{Q-1} 
p(k\theta +\varphi) 
\right|
\le 
 \sum_{\, 0<|m|\le N}
 |p_m|  \left|\frac{
 e^{2\pi i mQ\theta} -1}{  e^{2\pi i m\theta} -1}
\right|\le  \frac{\pi}{2} \frac{||Q\theta||}
{\min_{0<m \le N}  ||m\theta||} \sum_{\, 0<|m|\le N}
 |mp_m|.
$$
But by the best approximation property (\ref{be})   and  condition $ N< Q_{\nu+1}$ we have
$$
\min_{0<m \le N}  ||m\theta|| \ge ||Q_{\nu}\theta||
$$
and everything is proven.
\end{proof}

\subsection{ Sum with Kronecker sequence}\label{f4}

Here we prove that for any irrational $\theta$
for the function $f$ constructed in Subsection \ref{f}
one has 
 \begin{equation}\label{01}
 \liminf_{\nu\to \infty} \left| \sum_{k=0}^{Q_\nu-1} f(k\theta +\varphi) \right| = 0
 \end{equation}
and so 
(\ref{0}) is valid.

\vskip+0.3cm

We fix $\nu$ and consider $n=n(\nu)$ such that
\begin{equation}\label{e1}
\varepsilon_{n+2} <||Q_\nu\theta|| \le \varepsilon_{n+1}.
\end{equation}
Inequality (\ref{12}) together with right inequality from (\ref{e1})  and the upper bound $ \varepsilon_{n+1}.$ from (\ref{ee}) gives
$$
\frac{1}{2Q_{\nu+1}} <||Q_\nu\theta|| \le \varepsilon_{n+1}\le \frac{1}{2N_n}.
$$
So
\begin{equation}\label{nn}
N_n < Q_{\nu+1}.
\end{equation}
Then we consider  the minimal $l$ such that $ ||Q_l\theta|| \le \varepsilon_{n+3}$.
Then
\begin{equation}\label{e2}
||Q_l\theta|| \le \varepsilon_{n+3} < ||Q_{l-1}\theta|| .
\end{equation}
Now we are ready to finalise the proof of the theorem.
Write
 $$
f(x) = w(x) + p_n(x)+ h_{n+1}(x) + h_{n+2}(x) + g(x) , \,\,\, 
w(x) = f(x) - p_n(x),\,\,\, g(x) = \sum_{j=n+3}^\infty h_j(x)
$$
and divide the sum (\ref{sum}) into four  parts
\begin{equation}\label{s1}
\sum_{k=0}^{Q_l-1}
f(k\theta+\varphi) = S_1 + S_2 +S_3+S_4,
\end{equation}
where
$$
S_1 =  \sum_{k=0}^{Q_l-1}
w(k\theta+\varphi),
\,\,\,\,\,\,
S_2 =  \sum_{k=0}^{Q_l-1}
g_n(k\theta+\varphi)
$$
$$
S_3 =  \sum_{k=0}^{Q_l-1}
h_{n+1} (k\theta+\varphi)
+
 \sum_{k=0}^{Q_l-1}
h_{n+2} (k\theta+\varphi)
,
\,\,\,\,\,\,
S_4 =  \sum_{k=0}^{Q_l-1}
g(k\theta+\varphi).
$$
First of all we deal with sum $S_1$. As for the discrepancy of the sequence
$\Xi =\{ \{ k\theta\} , k = 0,...,Q_l-1\}$ we have the  bound
\begin{equation}\label{dva}
D_\Xi \le 2
\end{equation}
and $\int_0^1w(x) dx = 0$, Koksma's inequality (\ref{ko}) gives 
\begin{equation}\label{m1}
|S_1| \le 2 {\rm Var} [w] = 2\int_0^1|f_n'(x) - p_n'(x)|dx \le \frac{2}{n}  
\end{equation}
(we use (\ref{1})).
For the sum $S_2$  because of (\ref{nn}) we can use Lemma 2  for $p(x) = p_n(x) $ and $Q= Q_l$,  inequalities from (\ref{e1},\ref{e2})
and definition (\ref{2})   to get
\begin{equation}\label{m2}
|S_2 |\le
\frac{\pi}{2}\, 
M[p_n] \cdot \frac{||Q_l\theta||}
{||Q_{\nu}\theta||}\le
\frac{\pi}{2}\, 
M[p_n] \cdot \frac{\varepsilon_{n+3}}
{\varepsilon_{n+2}}\le \frac{\pi}{2(n+2)}
.
\end{equation}
For each sum form the  sum $S_3$   we use (\ref{dva}), Koksma's inequality (\ref{ko}) and bound
(\ref{var}) to get
\begin{equation}\label{m3}
|S_3|\le ({\rm Var}[h_{n+1}] + {\rm Var}[h_{n+2}] )D_\Xi\le
2 ({\rm Var}[h_{n+1}] + {\rm Var}[h_{n+2}] )\le 2\left(
\frac{4}{n+1}+\frac{4}{n+2}\right) \le \frac{16}{n}.
\end{equation}
For the sum $S_4$ we
observe that for $ x\not \in I = (x_{n+2},x^*)$
(the endpoints here are defined in (\ref{xx}))
we have $ g(x) = 0$.
We must tale into account that 
$$
|I| = x^* - x_{n+2} = \sum_{j=n+3} ^\infty \varepsilon_j \le 2\varepsilon_{n+3} < 2||Q_{l-1}\theta|| 
$$
by (\ref{ee}) and (\ref{e2}).
Now we can use  Lemma 1 and get
\begin{equation}\label{m4}
|S_4|\le  \sup_{x\in [0,1]}  \sum_{j=n+3}^\infty |h_j(x)|\le 
2\sup_{x\in [0,1]}  |h_{n+3}(x)| = \frac{2}{n+3}.
\end{equation}
Now we substitute (\ref{m1}) -- (\ref{m4}) into (\ref{s1}) and get
$$
\left|
\sum_{k=0}^{Q_l-1}
f(k\theta+\varphi) 
\right|
\le \frac{22}{n}.
$$
As $n(\nu) \to \infty$, because of $\lVert Q_\nu\theta \rVert \to 0$ when $\nu \to \infty$, everything is proven. $\Box$

\section{Proof of Theorem 2}

 Here we give a proof of Theorem 2.
 First of all, in Subsection \ref{l1}  we construct  a special function $g_\theta (x)$.
 Our argument is some kind of a generalization of the classical construction of a singular function. 
To study  properties of our sums we need to use the famous Three Distance Theorem.
We formulate all necessary facts about this theorem in Subsection \ref{l2}.
It would be convenient for us to use a quantitative version from \cite{be}.
A modern survey of the results related to the Three Distance Theorem one can find in \cite{mar}.
In Subsection \ref{l3} we study the behavior of the function $g_\theta (x)$
at the intervals where it is non-constant. In Subsection \ref{l4} we finalize the proof of Theorem 2.
 
\subsection {Construction of function $g_\theta (x)$}\label{l1}
Let irrational  $\theta$ be written as a continued fraction $\theta = [0;a_1,a_2,a_3,...]$ with $ a_1 = a_2 = 1$,
and  $Q_2 = 2$. In this section for every such $\theta$ we construct a special continuous function
$f(x)= f_\theta(x):[0,1]\to [0,1]$ such that 
on the interval $[0,\theta]$ function $f(x)$ is decreasing from 1 to 0, meanwhile
on the interval $[\theta, 1]$ function $f(x)$ is increasing from 0 to 1. In particular
${\rm Var}[f] =2 $.

Special properties of function $f(x)$ with respect to 
Kronecker sequence $\{k\theta\}, k = 1,2,3,...$ will be of importance in the next section.  Function $g_\theta(x)$ we define by
\begin{equation}\label{oogteta}
g_\theta (x) = f_\theta(x) - \int_0^1 f_\theta(x)dx.
\end{equation}
It is clear that $g_\theta (x)$ is continuous, 
on the interval $[0,\theta]$ it decreases, meanwhile
on the interval $[\theta, 1]$ it  increases,
${\rm Var}[g_\theta] =2 $
 and $g_\theta (x)$ satisfies equality (\ref{oozero}). 


\vskip+0.3cm

  First of all we  construct by induction all the values $f(x)$  when $x$ belongs to a special set $S$ which will be also constructed during our inductive process.

  The first step of the process corresponds to the value $ \nu=2$. At the first step we define just three values
\begin{equation}
    f(0) = f(1) = 1, f(\theta) = 0.
\end{equation}
Now let 
$$
f(0),f(\{\theta\}),...\, ,f(\{(Q_\nu-1)\theta\})  
$$
be already defined.
Let 
$$x_{\nu,0}, \, x_{\nu,1},\,  x_{\nu,2}\, , ... \, , \, x_{\nu,Q_\nu - 1}
$$
be permutation of the points 
$$
0, \, \{\theta\},\,  \{2\theta\}\, ,  ... \,, \, \{(Q_\nu - 1)\theta\}
$$ 
such that
\begin{equation}\label{oo2}
0 =x_{\nu,0} < x_{\nu,1} < x_{\nu,2} < ... < x_{\nu,Q_\nu - 1}
.
\end{equation}
We complete the sequence (\ref{oo2}) by the point $ x_{\nu, Q_\nu}= 1.$

Let $x_{\nu,Q_\nu} = 1$. Let us denote the segment $[x_{\nu,k-1}, x_{\nu, k}]$ by $I_{\nu, k}$, where $k = 1, 2, ..., Q_\nu$. We call $I_{\nu, k}$ a {\it constant segment} if $f(x_{\nu, k-1}) = f(x_{\nu, k})$, and a {\it jump segment} otherwise. At the first step $\nu = 2$ and we have only two jump segments, $[0, \theta]$ and $[\theta, 1]$.
However now we describe the rule of defining all the values of $f$ in a constant segment $I_{\nu,k} =
[x_{\nu, k-1},x_{\nu, k}].
$
The rule is as follows: if the values of $f(x_{\nu, k-1}) = f(x_{\nu, k})$ are defined at the previous step, we simply put
\begin{equation}\label{oo7}
f(x) = f(x_{\nu, k-1}) = f(x_{\nu, k})\,\,\,\text{for all}\,\,\, x\in I_{\nu,k} =
[x_{\nu, k-1},x_{\nu, k}].
\end{equation}
Now let us  define the set 
$$C_{\nu} := \{k: I_{\nu, k} \,\,\, \text{ is a  constant  segment}\}.$$
For the first step $\nu=2$ this set is empty, but it may be not empty in the next steps.
We should note that by inductive assumption 
the values 
$$
f(x_{\nu,0}), \, f(x_{\nu,1}), \,  f(x_{\nu,2}),  ...  \, , \, f(x_{\nu,Q_\nu - 1}),\, f(x_{\nu,Q_\nu }) 
$$
are already defined. So the rule described above defines $f(x) $ for all $x$ from the set 
$$S_\nu = \underset{k \in C_{\nu}}{\bigcup} I_{\nu, k}.$$
Now let 
proceed with the inductive step and define $f(x) $ at all points $x$ of the form
\begin{equation}\label{oo5}
\{Q_{\nu}\theta\}\,, \,
...\, ,\, \{(Q_{\nu+1} - 1)\theta\}.
\end{equation}
Some of the points (\ref{oo5}) fall into constant segments constructed up to  step $\nu$  and some of them fall into jump segments
constructed up to  step $\nu$.
If point $x$ of the form (\ref{oo5}) 
 falls into 
a constant segments $I_{\nu,k}$,
the value of 
$f(x)$ is already defined by induction hypothesis.
So it remains to determine $f(x)$ 
for those $x$ of the form (\ref{oo5})
  that fall into jump segments. 
  Let us fix a 
 jump segment $I_{\nu,k}$
 Either $I_{\nu,k}\subset [0,\theta]$ and $ x_{\nu,k}\le \theta$ or 
 $I_{\nu,k}\subset [\theta,1]$ and $ x_{\nu,k-1}\ge \theta$.
 Now we consider all the
 points $x$
 of the form (\ref{oo5}) which fall into  this segment $I_{\nu,k}$.
 For these points we define
\begin{equation}\label{oocase}
    \begin{cases}
    f(x) := \frac{1}{2}f(x_{\nu,k-1}) + \frac{1}{2}f(x_{\nu,k}), \ \text{if} \ x_{\nu,k-1} \ge \theta, \\
    f(x) := \frac{1}{3}f(x_{\nu,k-1}) + \frac{2}{3}f(x_{\nu,k}), \ \text{if} \ x_{\nu,k} \leq \theta.
    \end{cases}
\end{equation}
In such a way we define all the values 
$$f(\{Q_{\nu}\theta\})\, , ...\, ,\, f(\{(Q_{\nu+1} - 1)\theta\}).
$$
Now 
we can consider $(\nu+1)$-th collection of points 
$$
0 =x_{\nu+1,0} < x_{\nu+1,1} < x_{\nu+1,2} < ... < x_{\nu+1,Q_{\nu+1} - 1}<x_{\nu+1,Q_{\nu+1} } =1
$$
and the corresponding segments 
 $$
 I_{\nu+1,k} = [x_{\nu+1,k-1},x_{\nu+1,k}], \,\,\, k = 1, 2,...,Q_{\nu+1}.$$
 The value of $f(x)$ is defined at all endpoints of all segments $
 I_{\nu+1,k}$. Among these segments are 
    constant segments and jump segments. Some of the constant segments are subsets of other constant segments for previous steps and at the points of these segments the value of $f(x)$ is yet defined. Some of the constant segments of the $(\nu+1)$-step are new and we  define 
    the values of $f(x)$ for the points of these new constant segments by the rule (\ref{oo7}) described above.

 The inductive procedure described above define the values of $f(x)$ for all points x from the set    
      $$S = \bigcup\limits_{\nu = 2}^{\infty} S_\nu.$$
      At this moment we complete the inductive process. Now we must define the values of $f(x)$ for $x \not\in S$.
      We should note that any $x \not\in S$ falls into an interior of a jump segment. Moreover, ever such $x$ should belong to an intersection
      \begin{equation}\label{oointer}
      \bigcap_{r=1}^\infty I_{{\nu_r}, k_{\nu_r}},\,\,\,\, I_{2, k_2}\supset...\supset I_{\nu_r, k_{\nu_r}}\supset I_{\nu_{r+1}, k_{\nu_{r+1}}}\supset ... \, ,
    \end{equation}
      where each $ I_{\nu_r, k_{\nu_r}}$ is a jump segment 
      and
      $$
      2=\nu_1< ...<\nu_r<\nu_{r+1}<...\, .
      $$

      It is clear that $\mu(I_{\nu_r, k_{\nu_r}}) \rightarrow 0$
      for any choice of $ k_\nu$ when $\nu \rightarrow \infty$, because  Kronecker sequence is dense. So the intersection (\ref{oointer}) consists of just one point.

        For any segment $I = [\alpha, \beta]$ let us denote
$$\Delta_f(I) = |f(\alpha) - f(\beta)|.$$
Formulas (\ref{oocase}) show that for any different two nested jump segments $I_{\nu, k}\supset I_{\nu', k'} $
with $ \nu'>\nu$
we have
\begin{equation}\label{oonested}
\Delta_f(I_{\nu',k'}) \leq \frac{2}{3}\Delta_f(I_{\nu,k}).   
\end{equation}
So for any $x$ from the intersection (\ref{oointer}) we have 
$$|f(x_{\nu_r,k_{\nu_r}}) - f(x_{\nu_r,k_{\nu_r}-1})| \leq \bigg(\frac{2}{3}\bigg)^{r}.$$
We see that  $|f(x_{\nu_r,k_{\nu+r}}) - f(x_{\nu_r,k_{\nu_r}-1})| \rightarrow 0$ when $r \rightarrow \infty$. Also by construction $f$ decreases on the set   $S \cap [0, \theta]$ and increases on the set $S \cap [\theta, 1]$. It means that  for all $ x \in [0,1]$  both limits
$$ \lim_{S\ni y \to x-0} f(y)\,\,\,\text{and}\,\,\, \lim_{S\ni  y \to x+0} f(y)
$$
exist and are equal, that is
$$ \lim_{S\ni y \to x-0} f(y) = \lim_{S\ni y \to x+0} f(y)
= \lim_{S\ni y \to x } f(y).$$
Now for  $x \notin S$ we define 
$$f(x) = \lim_{S\ni y \to x} f(y).$$
Function $f$ defined in such  way 
will be   decreasing  from 1 to 0 on the segment $[0, \theta]$ and increasing from 0 to 1 on the segment  $[ \theta,1]$.  Moreover the limit conditions above show that $f(x)$ is continuous in every point $x \in [0,1]$.

  \subsection{Continued fractions and Three Distance Theorem}\label{l2}
  
For the analysis of  behavior of the sum 
\begin{equation}\label{nillsum}
\sum_{k=0}^{Q_{\nu}-1} g_\theta (\{k\theta\})
=
\sum_{k=0}^{Q_{\nu}-1} g_k
\end{equation}
 we use the following corollary of the famous {\it 
 Three Distance Theorem} in the very precise form given in \cite{be}
 (we use  Theorem 3 from \cite{be} in a particular case $ N= rq_k+q_{k-1}+s$ with
 $q_k = Q_{\nu-1}, q_{k-1} = Q_{\nu-2}$
and $ r= a_\nu-1,s = Q_{\nu-1}-1$; then $N_A =Q_\nu-Q_{\nu-1}, N_B = Q_{\nu-1}, N_C=0$, in notation of  \cite{be}).

\vskip+0.3cm
{\bf Lemma 3.}
{
Among segments $I_{\nu,k}$, where $k = 1,2,...,Q_\nu$ exactly $(Q_\nu - Q_{\nu - 1})$ segments of length 
$ 
||Q_{\nu-1}\theta||
$ and $Q_{\nu - 1}$ segments of length $||Q_{\nu-1}\theta|| + ||Q_\nu\theta||$. 
}

\vskip+0.3cm
 Also in the sequel we need the following well known equality
\begin{equation}\label{oo1000}
 \frac{||Q_{\nu-1}\theta||}{||Q_{\nu}\theta||} = \alpha_{\nu+1} = [ a_{\nu+1}; a_{\nu+2}, a_{\nu+3}, ...] .
 \end{equation}


\subsection{Function $g_\theta(x)$ on jump segments}\label{l3}




We continue with analysis of the properties of $g_\theta(x)$.
 Consider points $x_{\nu,k}$ and  segments $I_{\nu,k}$ defined in the previous section.
Let us denote $g_\theta(x_{\nu,k}) = g_k$. It is easy to see that if $I_{\nu,k}$ is a constant segment, then
\begin{equation}
    \int\limits_{I_{\nu,k}} g_\theta(x)dx = g_{k-1}\mu(I_{\nu,k}).
\end{equation}
Define the sets of jump segments
$$J_{\text{left}} := \{k : I_{\nu,k+1} \,\,\text{ is  a  jump  segment and}\,\,\, x_{\nu, k+1} \leq \theta\}$$
and
$$J_{\text{right}} := \{k : I_{\nu,k+1} \,\,\text{ is  a  jump  segment and}\,\,\,x_{\nu, k} \geq \theta\}.$$ 
Now let us estimate the integrals over the jump segments.
We distinguish the cases when the jump segment is form the set $J_{\text{left}}$ and 
when the jump segment is form the set $J_{\text{right}}$.

\vskip+0.3cm
{\bf 1.  First  case}, $k -1\in J_{\text{left}}$. 

In this case
  by construction of $g_\theta$ on the segment $I_{\nu,k}$  function $g_\theta(x)$ decreases monotonically.
  Consider the set 
$$
  \Omega_+ = 
  \left\{x \in I_{\nu,k}: g_\theta (x) >  G_k^{\text{left}}\right\}
  ,$$
  where
  \begin{equation}\label{oogk}
  G_k^{\text{left}}=
  \frac{1}{3}g_{k-1} + \frac{2}{3}g_k,\,\,\,\, g_{k-1}\geq G_k^{\text{left}}\geq g_k.
  \end{equation}
  We prove that 
\begin{equation}\label{oomu1}
\mu (\Omega_+) 
< 2 ||Q_{\nu}\theta||.
\end{equation}
Indeed, $I_{\nu,k}=[x_{\nu,k-1},x_{\nu,k}]$ is a jump segment. Consider the smallest $\nu'>\nu$ such that a point of the form $x_{\nu',k'}$  
falls into $I_{\nu,k}$ for the first time. Let this point $x_{\nu',k'}$ be just the minimal point which falls into $I_{\nu,k}$.
Then by our construction
$$f_\theta (x_{\nu',k'}) = \frac{1}{3}f_\theta(x_{\nu,k-1}) + \frac{2}{3}f_\theta (x_{\nu,k}).$$
By the definition (\ref{oogteta}) we see that also
$$g_\theta (x_{\nu',k'}) = \frac{1}{3}g_\theta(x_{\nu,k-1}) + \frac{2}{3}g_\theta (x_{\nu,k})=
G_k^{\text{left}}.$$
So by monotonicity
$$
\Omega_+ =
[x_{\nu,k-1}, x_{\nu',k'}).$$
By Lemma 3 (or by Three Distance Theorem) the length of the interval $[x_{\nu,k-1}, x_{\nu',k'}]$ is
either  $||Q_{\nu'-1}\theta||$ or  $||Q_{\nu'-1}\theta||+||Q_{\nu'}\theta||$. But $\nu'-1\geq \nu$ and so 
$$
\max (
||Q_{\nu'-1}\theta||, ||Q_{\nu'-1}\theta||+||Q_{\nu'}\theta||) < 2 ||Q_{\nu'-1}\theta||\leq 2||Q_{\nu}\theta||.
$$
Inequality (\ref{oomu1}) is proven.

\vskip+0.3cm

Now for the integral over $I_{\nu,k}$ we can write the equality
 $
    \int\limits_{I_{\nu,k}}   =   \int \limits_{I_{\nu,k}\setminus \Omega_+}+ \int\limits_{\Omega_+} .
   $
When $x\in I_{\nu,k}\setminus \Omega_+ $ we have the bound  $g(x)\leq\frac{1}{3}g_{k-1} + \frac{2}{3}g_k$
 , meanwhile when 
$x \in I_{\nu,k}\setminus \Omega$ we have a trivial bound $ g(x)\leq \max_{I_{\nu,k}} g(x) = g_{k-1}.$
So, taking into account (\ref{oomu1}) and the inequality $g_{k-1}- G_k^{\text{left}}\ge 0$ which follows from (\ref{oogk}), we get the upper bound
   
$$
   \int\limits_{I_{\nu,k}} g_\theta(x)dx
\leq \int \limits_{I_{\nu,k}\setminus \Omega_+}G_k^{\text{left}} dx+ \int\limits_{\Omega_+}g_{k-1}dx =
G_k^{\text{left}} \cdot (\mu ( I_{\nu,k}) - \mu(\Omega_+)) + g_{k-1} \cdot \mu(\Omega_+) =
$$
     \begin{equation}\label{oo21}
    =
    G_k^{\text{left}} \cdot \mu ( I_{\nu,k}) +(  g_{k-1}- G_k^{\text{left}})\cdot \mu(\Omega_+){\leq}
    \bigg(G_k^{\text{left}} + \frac{4}{3}(g_{k-1} - g_k)\frac{||Q_\nu\theta||}{\mu(I_{\nu,k})}\bigg)\mu(I_{\nu,k}).
\end{equation}
Similarly to (\ref{oomu1}) from Lemma 3    we deduce inequality 
$$ \mu (\Omega_-) 
< 2 ||Q_{\nu}\theta||\,\,\,
\text{where}\,\,\,\Omega_- \left\{x \in I_{\nu,k}: g_\theta (x) < G_k^{\text{left}}\right\} < 2||Q_\nu\theta||,$$
which leads  by (\ref{oomu1}) and  (\ref{oogk}) to 
$$ 
    \int\limits_{I_{\nu,k}} g_\theta(x)dx \geq
    G_k^{\text{left}} \cdot (\mu ( I_{\nu,k}) - \mu(\Omega_-)) + g_{k} \cdot \mu(\Omega) =
    G_k^{\text{left}}\cdot \mu ( I_{\nu,k})- (G_k^{\text{left}} -g_k)\cdot\mu(\Omega_-)\geq
    $$
    \begin{equation}\label{oo22}
    \geq  G_k^{\text{left}}\cdot \mu ( I_{\nu,k})-  2(G_k^{\text{left}} -g_k)\cdot||Q_\nu\theta||
    = \bigg(G_k^{\text{left}} - \frac{2}{3}(g_{k-1} - g_k)\frac{||Q_\nu\theta||}{\mu(I_{\nu,k})}\bigg)\mu(I_{\nu,k}).
\end{equation}
Finally we collect (\ref{oo21}) and (\ref{oo22}) together and get for any $I_{\nu,k} \in J_{\text{left}} $ a double inequality
\begin{equation}\label{oo33}
\bigg(G_k^{\text{left}} - \frac{2}{3}(g_{k-1} - g_k)\frac{||Q_\nu\theta||}{\mu(I_{\nu,k})}\bigg)\mu(I_{\nu,k})\leq
\int\limits_{I_{\nu,k}} g_\theta(x)dx
\leq
\bigg(G_k^{\text{left}} + \frac{4}{3}(g_{k-1} - g_k)\frac{||Q_\nu\theta||}{\mu(I_{\nu,k})}\bigg)\mu(I_{\nu,k})
\end{equation}

\vskip+0.3cm

{\bf 2.  Second  case}, $k -1 \in J_{\text{right}}$. 

 Recall that  $g_\theta(x)$ increases monotonically on every $I_{\nu,k}\subset [\theta,1]$.
 For any jump segment $I_{\nu,k}$ with $k \in J_{\text{right}}$ by means of  values $g_k = g_\theta(x_{\nu,k})$ we define
 \begin{equation}\label{oorightw}
 G_k^{\text{right}} = \frac{1}{2}g_{k-1} + \frac{1}{2}g_k,\,\,\,\,  g_{k-1} \leq  G_k^{\text{right}}\leq  g_k.
 \end{equation}
 Then we repeat all the argument of the first case word by word and obtain a double inequality 
\begin{equation}\label{ooright1}
    \bigg( G_k^{\text{right}}- (g_{k} - g_{k-1})\frac{||Q_\nu\theta||}{\mu(I_{\nu,k})}\bigg)\mu(I_{\nu,k}) { \leq} \int\limits_{I_{\nu,k}} g_\theta(x)dx { \leq} \bigg(G_k^{\text{right}} + (g_{k} - g_{k-1})\frac{||Q_\nu\theta||}{\mu(I_{\nu,k})}\bigg) \mu(I_{\nu,k})
\end{equation}

\subsection{Proof of Theorem 2: the end}\label{l4}

For every  $k = 0, 2, ..., Q_{\nu}-1$ we consider the values
$$\widetilde g_{k}= \frac{1}{\mu(I_{\nu,k+1})}\int\limits_{I_{\nu,k+1}} g_\theta(x)dx.$$
The values $\widetilde g_{k}$ admit the following properties.
First of all, from the definition of $ g_{k}$ and $\widetilde g_{k}$ we see that 
\begin{equation}\label{gkconstant}
\widetilde g_k = g_k, \,\,\text{if} \,\, I_{\nu,k+1} \,\,\text{ is  a  constant  segment}. 
\end{equation}
Then, we should note that by Lemma 3 in any case 
$$ 
\mu(I_{\nu,k}) \ge ||Q_{\nu-1}\theta||.
$$
So from (\ref{oo33}) we see that 
\begin{equation}\label{gkleft2}
| \widetilde g_k - G_{k+1}^{\text{left}}| < \frac{4}{3} \,\frac{||Q_\nu\theta||}{||Q_{\nu-1}\theta||}\, |g_k - g_{k+1}|,
\,\,\text{if}\,\,I_{\nu,k+1} \,\,\text{ is  a  jump  segment and}\,\,
  k \in J_{\text{left}}. 
\end{equation}
Similarly, from (\ref{ooright1}) we see that 
\begin{equation}\label{gkright2}
|\widetilde g_k - G_{k+1}^{\text{right}}| < \frac{||Q_\nu\theta||}{||Q_{\nu-1}\theta||}\, |g_k - g_{k+1}|,
\,\,\text{if}\,\,I_{\nu,k+1} \,\,\text{ is  a  jump  segment and}\,\,
 k \in J_{\text{right}}. 
\end{equation}

From (\ref{oozero}) it follows that 
\begin{equation}\label{oosum}
    \sum_{k=0}^{Q_\nu-1}\widetilde g_k\mu(I_{\nu, k+1}) = \int_0^1 g_\theta(x)dx = 0
\end{equation}
We divide this quantity  by $||Q_{\nu-1}\theta||$. By Lemma 3 
in the sum  (\ref{oosum}) after the division 
will consist of 
$(Q_{\nu} - Q_{\nu-1})$
summands of the form $\widetilde g_k$ and
  and $Q_{\nu-1}$ 
  summands of the form 
  $\left( 1 + \frac{||Q_\nu\theta||}{||Q_{\nu-1}\theta||}\right) \widetilde g_k$,
  that is
  $$
  0 =
  \frac{1}{||Q_{\nu-1}\theta||}
   \sum_{k=0}^{Q_\nu-1}\widetilde g_k\mu(I_{\nu, k+1})=
    \Sigma^{(1)}\widetilde g_k +
   \Sigma^{(2)} \left( 1+\frac{||Q_\nu\theta||}{||Q_{\nu-1}\theta||}\right)\widetilde g_k,
  $$
  where in sum $ \Sigma^{(1)}$ we collect the summands with $\mu(I_{\nu,k+1}) = ||Q_{\nu-1}\theta||$ and in 
  sum $ \Sigma^{(2)}$ we collect the summands with $\mu(I_{\nu,k+1}) = ||Q_{\nu-1}\theta||+ ||Q_{\nu}\theta||.$
    Thus, because of $|\widetilde g_k| \leq 1$ we get
\begin{equation}\label{suumm}
    \bigg|\sum_{k=0}^{Q_\nu-1}\widetilde g_k \bigg|=
    \bigg|
     \Sigma^{(2)} \frac{||Q_\nu\theta||}{||Q_{\nu-1}\theta||}\widetilde g_k \bigg| \leq Q_{\nu - 1} \frac{||Q_\nu\theta||}{||Q_{\nu-1}\theta||}.
\end{equation}
Now combining 
(\ref{gkconstant},\ref{gkleft2},\ref{gkright2})
we get
\begin{equation}\label{boound}
\left|
\sum_{k=0}^{Q_\nu-1}\widetilde g_k -
\sum_{k: I(\nu,{k+1}) -\,\text{constant segment}} g_k
-
\sum_{k\in J_{\text{left}}} G_{k+1}^{\text{left}} -
\sum_{k\in J_{\text{right}}} G_{k+1}^{\text{right}}
\right| \le 
\frac{4}{3}\,\frac{||Q_\nu\theta||}{||Q_{\nu-1}\theta||}\sum_{k=0}^{Q_\nu-1}|g_k - g_{k+1}|.
\end{equation}
By the definitions (\ref{oogk},\ref{oorightw}) of $G_{k}^{\text{left}}$ and $G_{k}^{\text{right}}$ the left hand side here may be  rewritten as
$$
\sum_{k=0}^{Q_\nu-1}\widetilde g_k -
\sum_{k: I(\nu,{k+1}) -\,\text{constant segment}} g_k
-
\sum_{k\in J_{\text{left}}} G_{k+1}^{\text{left}} -
\sum_{k\in J_{\text{right}}} G_{k+1}^{\text{right}}=
$$
$$
=
\sum_{k=0}^{Q_\nu-1}\widetilde g_k -
\sum_{k: I(\nu,{k+1}) -\,\text{constant segment}} g_k
-
\sum_{k\in J_{\text{left}}} \left( \frac{1}{3} g_k + \frac{2}{3}g_{k+1}\right) -
\sum_{k\in J_{\text{right}}}
\left( \frac{1}{2} g_k + \frac{1}{2}g_{k+1}\right)=
$$
$$
=
\sum_{k=0}^{Q_\nu-1}\widetilde g_k -
\sum_{k=0}^{Q_\nu-1} g_k {+} \frac{2}{3}\sum_{k\in J_{\text{left}}}(g_k - g_{k+1}) {+} \frac{1}{2}\sum_{k\in J_{\text{right}}}(g_k - g_{k+1}).
$$
Taking into account inequality
$$
\sum_{k=0}^{Q_\nu-1}|g_k - g_{k+1}|\leq {\rm Var} [g_\theta]= 2
$$
we transform (\ref{boound}) into
\begin{equation} \label{booound}
    \bigg|\sum_{k=0}^{Q_\nu-1}g_k - \sum_{k=0}^{Q_\nu-1}\widetilde g_k - \frac{2}{3}\sum_{k\in J_{\text{left}}}(g_k - g_{k+1}) - \frac{1}{2}\sum_{k\in J_{\text{right}}}(g_k - g_{k+1})\bigg| <   \frac{8}{3}\, \frac{||Q_\nu\theta||}{||Q_{\nu-1}\theta||}.
\end{equation}
Now we should  note that
$$\sum_{k\in J_{\text{left}}}(g_k - g_{k+1}) = -\sum_{k\in J_{\text{right}}}(g_k - g_{k+1}) = g_{\theta}(0) - g_{\theta}(\theta) = 1,$$
and instead of (\ref{booound}) we can write
\begin{equation} \label{bnd}
    \bigg|\sum_{k=0}^{Q_\nu-1}g_k - \sum_{k=0}^{Q_\nu-1}\widetilde g_k - \frac{1}{6}\bigg| <   \frac{8}{3}\, \frac{||Q_\nu\theta||}{||Q_{\nu-1}\theta||}.
\end{equation}
Combining (\ref{bnd}) with  (\ref{suumm})  and (\ref{oo1000}) we get
$$
    \bigg|\sum_{k=0}^{Q_\nu-1}g_k - \frac{1}{6}\bigg| < 
    \bigg|\sum_{k=0}^{Q_\nu-1}\widetilde g_k\bigg| +  \frac{8}{3}\frac{D_\nu}{D_{\nu-1}}
    \leq \left(Q_{\nu-1} + \frac{8}{3}\right)\frac{||Q_\nu\theta||}{||Q_{\nu-1}\theta||} <4 \,  \frac{Q_{\nu-1} ||Q_{\nu}\theta||}{||Q_{\nu-1}\theta||}= \frac{4 Q_{\nu-1}}{\alpha_{\nu+1}}.
$$    
But but the assumption (\ref{oo1}) of Theorem 1?
$\frac{4 Q_{\nu-1}}{\alpha_{\nu+1}}\to 0, \nu \to \infty.$
So
$$
\lim_{\nu \to \infty}
\sum_{k=0}^{Q_{\nu}-1} g_\theta (\{k\theta\})
=
\lim_{\nu \to \infty} \sum_{k=0}^{Q_\nu-1}g_k = \frac{1}{6} > 0,$$
and Theorem 2 is proven.$\Box$

\section{Proof of Proposition 1}
 We may restrict ourselves to the case $ \varphi = 0$.
As 
 $f$ 
 is a function of bounded variation, condition (\ref{ooBV}) ensures  existence of $\gamma\in \mathbb{R}$ and a subsequence 
 $\{\nu_m\}$ , $ m=1,2,3,... $ such that simultaneously  condition (\ref{oo1o}) holds 
 and asymptotic equality  
 \begin{equation}\label{oooo1}
 \lim_{m\to \infty}\sum_{k=0}^{Q_{\nu_m}-1} f(k\theta) = \gamma
 \end{equation}
 is valid.
 Now we choose a subsequence $\nu_{n(m)}$ in such a way that 
 $$
 Q(m) =  Q_{\nu_{n(m)}}- Q_{\nu_{m}} \to \infty,\,\, m\to \infty.
  $$
  Note that 
  by the definition of the subsequence  $\{\nu_m\}$, we have
  \begin{equation}\label{ooo5}
 \sum_{k=0}^{Q_{\nu_{n(m)}-1}} f(k\theta) -
   \sum_{k=0}^{Q_{\nu_m}-1} f(k\theta) 
  \to \gamma - \gamma = 0,\,\,\, m \to \infty.
  \end{equation}
  Moreover,
 $$ 
  \sum_{k=Q(m)}^{Q_{\nu_{n(m)}}-1} f(k\theta) 
  =
  \sum_{k=0}^{Q_{\nu_{m}}-1} f(k\theta + Q_{\nu_{n(m)}}\theta - Q_{\nu_{m}}\theta) 
.
 $$ 
 But now  by the definition of the modulus of continuity $\omega(\cdot)$ and (\ref{oo1o}) we get 
 $$
\left|
   \sum_{k=Q(m)}^{Q_{\nu_{n(m)}}-1} f(k\theta) 
   -
   \sum_{k=0}^{Q_{\nu_{m}}-1} f(k\theta ) 
\right|
   =
  \left|
  \sum_{k=0}^{Q_{\nu_{m}}-1} f(k\theta + Q_{\nu_{n(m)}}\theta - Q_{\nu_{m}}\theta)
  -
    \sum_{k=0}^{Q_{\nu_{m}}-1} f(k\theta ) 
  \right|\leq
  $$
    \begin{equation}\label{ooo2}
  \leq Q_{\nu_{m}}
  (\omega (||Q_{\nu_{n(m)}}\theta||)
  +
  \omega (||Q_{\nu_{m}}\theta||))
  \leq
  2Q_{\nu_{m}} \omega (||Q_{\nu_{m}}\theta||)\to 0,\,\,\,
  m \to \infty  
.
  \end{equation}
  Taking into account (\ref{ooo5}) and (\ref{ooo2}) we conclude that 
  
  $$
  \sum_{k=0}^{Q(m)-1} f(k\theta) =
\sum_{k=0}^{Q_{\nu_{n(m)}}-1} f(k\theta)-
\sum_{k=Q(m)}^{Q_{\nu_{n(m)}}-1} f(k\theta)\to 0,\,\,\,\,
m\to \infty,
$$
 and we get (\ref{ooo1}).
 Proposition 1 is proven. $\Box$

\section{Proof of Proposition 2}

Similarly to the proof of Proposition 1, there exist
$\gamma\in \mathbb{R}$ and a subsequence 
 $\{\nu_m\}$ , $ m=1,2,3,... $ such that  
 $$
  \lim_{m\to \infty}\sum_{k=0}^{Q_{\nu_m}} f(k\theta) = \gamma
  .
  $$
 Let $\varepsilon > 0$ and $m \in \mathbb N$ be such that
\begin{equation}\label{ooooo7}
\left|\sum_{k=0}^{Q_{\nu_m}}f(k\theta) - \gamma \right| < \varepsilon    
\end{equation}
Since $f$ is   uniformly continuous,  there exists $\delta =\delta(\varepsilon, m)  > 0$ such that 
\begin{equation}\label{ooooo1}
|f(x) - f(y)| < {\varepsilon}{Q_{\nu_m}^{-1} }\,\,\,
\text{
for all }\,\,\, x,y \,\,\,
\text{with}\,\,\,
  |x-y| < \delta.
\end{equation}
  Since $\lVert Q_{\nu_m}\theta \rVert \to 0$ when $m \to \infty$ there exists $n = n(m)$ such that  simultaneousely we have 
\begin{equation}\label{ooooo1o}
\lVert Q_{\nu_{n(m)}}\theta \rVert < \delta
\end{equation}
and
\begin{equation}\label{ooooo6}
\left|\sum_{k=0}^{Q_{\nu_{n(m)}}}f(k\theta) - \gamma \right| < \varepsilon  .  
\end{equation}
Again we put
$
Q(m) =Q_{\nu_{n(m)}} - Q_{\nu_m} $.
We may assume that $ Q(m) \to \infty$ as $m\to\infty$.
From our assumption (\ref{even}) we have
\begin{equation}\label{ooooo3}
\sum_{k=Q(m)}^{Q_{\nu_{m}}}f(k\theta)=
\sum_{k=0}^{Q_{\nu_{m}}}f(Q_{\nu_{n(m)}}\theta -k\theta)=
\sum_{k=0}^{Q_{\nu_{m}}}f(k\theta- Q_{\nu_{n(m)}}\theta).
\end{equation}
By (\ref{ooooo1})  and  (\ref{ooooo1o}) we have
\begin{equation}\label{ooooo9}
    \left |\sum_{k=0}^{Q_{\nu_{m}}}f(k\theta - Q_{\nu_{n(m)}}\theta) - \sum_{k=0}^{Q_{\nu_m}}f(k\theta  )\right |\leq (Q_{\nu_m} +1)\cdot
    \max_x   \left| f(x - Q_{\nu_{n(m)}}\theta) - f(x)\right|
  \leq  2\varepsilon.
\end{equation}
So, applying (\ref{ooooo7},\ref{ooooo6},\ref{ooooo3},\ref{ooooo9}) we get
$$ 
    \left |\sum_{k=0}^{Q(m)-1}f(k\theta)\right| = \left|\sum_{k=0}^{Q_{\nu_{n(m)}}}f(k\theta) \ - \sum_{k = Q(m)}^{Q_{\nu_{n(m)}}}f(k\theta)\right| < 4\varepsilon.
$$
This proves Proposition 2.$\Box$

\section{Acknowledgements}
I express my gratitude to R.F. Tichy for setting the question and to N. Moshchevitin for his close attention to the work and valuable comments. The work was supported by a grant from the Foundation for the Advancement of Theoretical Physics and Mathematics "BASIS" and the Moscow Center for Fundamental and Applied Mathematics.

\end{document}